\documentclass[
final
 , nomarks
]{dmtcs-episciences}

\usepackage[utf8]{inputenc}
\usepackage{subfigure}
\usepackage{amsmath,amsthm,amssymb}
\newcommand{\seq}[1]{\ensuremath{\left( #1 \right)}}
\newtheorem{theorem}{Theorem}

\usepackage[round]{natbib}

\author{Brian Cloteaux\affiliationmark{1}}
\title[A Sufficient Condition for Graphic Sequences]{
A Sufficient Condition for Graphic Sequences with Given Largest and Smallest
Entries, Length, and Sum
\footnote{Official contribution of the National Institute of Standards and Technology; not subject to copyright in the United States.}
}

\affiliation{
  National Institute of Standards and Technology, Applied and
  Computational Mathematics Division, USA}

\received{2016-10-19}
\revised{2017-3-9, 2018-1-22, 2018-6-8, 2018-6-13}
\accepted{2018-06-14}
\begin{document}
\publicationdetails{20}{2018}{1}{25}{4600}
\maketitle
\begin{abstract}
We give a sufficient condition for a nonnegative integer list to be
graphic based on its largest and smallest elements, length, and sum.  This
bound generalizes a result of Zverovich and Zverovich.
\end{abstract}

We denote a finite list of nonnegative integers,
$\alpha=\seq{\alpha_1, ..., \alpha_n}$ where
$\alpha_1 \geq ... \geq \alpha_n$ as a {\it degree sequence}.
A degree sequence is said
to be {\it graphic} if it corresponds to the set of edge adjacency values of the
nodes for some simple graph.  
In 1960, Erd\H{o}s and Gallai gave a complete
characterization of the set of graphic degree sequences.
\begin{theorem}[\cite{Erdos:1960}]
\label{thm:erdos-gallai}
Let $\alpha=\seq{\alpha_1, ..., \alpha_n}$ be a nonincreasing degree sequence.
Then the degree sequence $\alpha$ is graphic if and only if the
sum of $\alpha$ is even and for each integer $k$ where $1 \leq
k \leq n$,
\begin{equation}
\sum^{k}_{i=1} \alpha_i \leq k(k-1) + \sum^{n}_{i=k+1} \min \lbrace k, \alpha_i
\rbrace.
\end{equation}
\end{theorem}

It was later shown by I. Zverovich and E. Zverovich that 
some degree sequences with bounded largest and smallest elements can be
verified to be graphic based on their length. 

\begin{theorem} [\cite{Zverovich:1992}, Theorem 6]
\label{thm:zz}
Let $\alpha=\seq{\alpha_1, ..., \alpha_n}$ be a nonincreasing degree sequence of
positive integers with even sum.  If
\begin{equation}
n \geq \frac{(\alpha_1 + \alpha_n + 1)^2}{4\alpha_n},
\end{equation}
then $\alpha$ is graphic.
\end{theorem}

There have been several extensions to the original result of Zverovich
and Zverovich.
These results have been either specializations, such as
when the gaps between consecutive integers in the sequence are bounded 
(\cite{Barrus:2012}),
or by sharpening the bound for specific sequences
(\cite{Cairns:2015,Cairns:2016}).
In this article, we show the following 
generalization of the Zverovich and Zverovich result by taking into account
the sequence sum.

\begin{theorem}
\label{thm:edge_bound}
Let $\alpha=\seq{\alpha_1, ..., \alpha_n}$ be an integer sequence 
such that $n-1 \geq \alpha_1 \geq ... \geq \alpha_n \geq 0$  and with even
sum $s = \sum_{i=1}^{n} \alpha_i$ where $n\alpha_1 > s > n\alpha_n$.
If
\begin{equation} \label{eqn:edge_bound}
(\alpha_1 - \alpha_n) \Big( \frac{n-\alpha_1-1}{n\alpha_1 - s} +
\frac{\alpha_n}{s - n\alpha_n} \Big) \geq 1,
\end{equation} 
then $\alpha$ is graphic.
\end{theorem}

\section{Proof of Theorem \ref{thm:edge_bound}}
Before providing a proof of this result, we state some needed
definitions and theorems.
The {\it complement} of a degree sequence
$\alpha = \left( \alpha_1, \alpha_2, ..., \alpha_n \right)$
is the sequence $\bar{\alpha} = \left( n-\alpha_n-1, ..., n-\alpha_1-1 \right)$.
It is straightforward to see that $\alpha$ is graphic if and only if $\bar{\alpha}$ 
is graphic.  It is important to note that Equation \eqref{eqn:edge_bound} is
invariant under complementation, meaning that 
\begin{equation} 
(\alpha_1 - \alpha_n) \Big( \frac{n-\alpha_1-1}{n\alpha_1 - s} +
\frac{\alpha_n}{s - n\alpha_n}
\Big)  = 
(\bar{\alpha_1} - \bar{\alpha_n}) \Big(
\frac{n-\bar{\alpha_1}-1}{n\bar{\alpha_1} - \bar{s}} +
\frac{\bar{\alpha_n}}{\bar{s} - n\bar{\alpha_n}}
\Big),
\end{equation} 
where $\bar{s} = \sum_{i=1}^{n} \bar{\alpha}_i$.

A second definition is that a sequence $\alpha$ {\it majorizes} a second
sequence $\beta = \seq{\beta_1, ..., \beta_n}$, where both sequence have the
same length and sum, if and only if
\begin{equation}
\label{eqn:major_def}
\sum_{i=1}^{k} \alpha_i \geq \sum_{i=1}^k \beta_i,
\end{equation}
for each $k$ from $1$ to $n$.

Majorization is a partial order over the set of degree sequences with
identical sum.  Our use of majorization stems from the following result. 
\begin{theorem}[\cite{Ruch:1979}, Theorem 1]
\label{thm:graphicality}
If the degree sequence $\alpha$ is graphic and $\alpha$ majorizes $\beta$,
then $\beta$ is graphic.
\end{theorem}

In addition, we use the following two results that reduce the number of
Erd\H{o}s-Gallai inequalities needed to verify that a sequence is
graphic.
These are found in \cite{Tripathi:2003}.
\begin{theorem}[\cite{Tripathi:2003}, Theorem]
\label{thm:run_check}
For the degree sequence $\alpha=(\alpha_1,...,\alpha_n)$, define the 
set of indices $\mathcal{I}$ such that $i \in \mathcal{I}$ if $\alpha_i > \alpha_{i+1}$. 
For Theorem \ref{thm:erdos-gallai}, it suffices to check only the
inequalities for the indices in $\mathcal{I}$.
\end{theorem}

\begin{theorem}[\cite{Tripathi:2003}, Lemma]
\label{thm:dur_index}
For the degree sequence $\alpha=(\alpha_1,...,\alpha_n)$, define the 
set of indices $\mathcal{I}$ where $i \in \mathcal{I}$ if $\alpha_i \geq i-1$. 
For Theorem \ref{thm:erdos-gallai}, it suffices to check only the
inequalities for the indices in $\mathcal{I}$.
\end{theorem}

We are now ready to prove Theorem \ref{thm:edge_bound}.

\begin{proof}[of Main Theorem]
We begin by defining
the set $\mathcal{D}$ as the set of all  degree sequences
that satisfy the hypotheses conditions for $\alpha_1$, $\alpha_n$,
$n$, and $s$.  We then define the sequence $\alpha' \in \mathcal{D}$ as 
\begin{equation}
\label{eqn:max_maj}
\alpha' = \left( \underbrace{\alpha_1, ..., \alpha_1}_{p},\gamma,
	\underbrace{\alpha_n, ..., \alpha_n}_{n-p-1} \right),
\end{equation}
where $\alpha_1 > \gamma \geq \alpha_n$ and $s = p\alpha_1 + \gamma + (n-p-1)\alpha_n$.
From the condition that $n\alpha_1 > s > n\alpha_n$, it follows that the
sequence $\alpha'$ cannot be regular, i.e., $\alpha_1 > \alpha_n$.
To show that the sequence $\alpha'$ exists, we construct it by first defining the
value of $\gamma$ as $\gamma = \tau + \alpha_n$ where $\tau$ is the solution to the
congruence equation 
\begin{equation}
	\tau \equiv s - n\alpha_n  \pmod{(\alpha_1 - \alpha_n)}, 
\end{equation}
that satisfies $0 \leq \tau < \alpha_1 - \alpha_n$. Since $\alpha_1 >
\alpha_n$, then $\tau$ is well-defined. 
From definition of $\tau$, the value of $\gamma$ is 
$\alpha_n \leq \gamma  < \alpha_1$. It follows from the congruence equation that
$s - n\alpha_n - \tau = p(\alpha_1 - \alpha_n)$ for some
integer $p$.  
To show that the sequence $\alpha'$ can be constructed, we simply need to
show that $0 < p < n$.  First, we rewrite the value for $p$ as
\begin{equation}
	p = \frac{s - n\alpha_n - \tau}{\alpha_1 - \alpha_n}.
\end{equation}
From $0 \leq \tau < \alpha_1-\alpha_n$, $\alpha_n \leq ... \leq \alpha_1$, and
$\sum_{i=1}^{n} \alpha_i = s < n\alpha_1$, then
\begin{equation} \label{eqn:p_inequality}
0 < \frac{\alpha_1-\alpha_n-\tau}{\alpha_1-\alpha_n} \leq 
\frac{\alpha_1 + (n-1)\alpha_n - n\alpha_n - \tau}{\alpha_1-\alpha_n} \leq
\frac{s - n\alpha_n - \tau}{\alpha_1-\alpha_n} <
\frac{n\alpha_1 - n\alpha_n - \tau}{\alpha_1-\alpha_n} \leq
n,
\end{equation}
implying $0 < p < n$ and establishing that the sequence $\alpha'$ does
uniquely exists.

We observe that the sequence
$\alpha'$ majorizes all the other sequences in $\mathcal{D}$.
If there would exist a sequence
$\beta \in \mathcal{D}$ where $\alpha'$ does not majorize $\beta$ then
there would exist an index $k$ where
$\sum_{i=1}^{k} \beta_i > \sum_{i=1}^k \alpha'_i$. The first index
where this could
occur is at $k = p+1$, since the first $p$ values of $\alpha'$ are the maximum
possible value $\alpha_1$.
If $\sum_{i=1}^{p+1} \beta_i > \sum_{i=1}^{p+1} \alpha'_i$, then
$\sum_{i=p+2}^{n} \beta_i <\sum_{i=p+2}^{n} \alpha' = (n+p-1)\alpha_n$,
forcing a value in the sequence $\beta$ to be less than
$\alpha_n$ and violating its membership in the set $\mathcal{D}$. There is
a similar argument for remainder of the indices, establishing that $\alpha'$
must majorize all the members of $\mathcal{D}$.
It  follows from Theorem \ref{thm:graphicality} that if $\alpha'$ is graphic
then every sequence in $\mathcal{D}$ must also be graphic. Thus we need only
to consider whether the sequence $\alpha'$ is graphic in order to prove the
result.

We also make the assumption concerning $\alpha'$ that $\alpha_1 \geq
n-\alpha_n-1 = \bar{\alpha_1}$.
This assumption holds in general because if $\alpha_1 < \bar{\alpha_1}$
then we simply use the complement of $\alpha'$ for our argument since
Equation \eqref{eqn:edge_bound} is invariant under complementation.

Now let us examine the  Erd\H{o}s-Gallai (EG) inequalities for $\alpha'$.
We apply Theorem \ref{thm:run_check},
to show that if the EG inequalities hold at $k=p$ and
$k=p+1$ in order to prove that the sequence $\alpha'$ is graphic.

First consider the case when $p < \alpha_n$.  The resulting EG
inequality for $k=p$  is
\begin{equation}
	p\alpha_1 \leq p(p-1) + p(n-p) = p(n-1),
\end{equation}
which is false only when $\alpha_1 > n-1$, but we assumed that $\alpha_1 \leq n-1$
ensuring that this inequality always holds.
There is an identical argument for the case when $k=p+1$ establishing that the
EG inequalities always hold when $p < \alpha_n$.

Now consider when the case when $p \geq \alpha_n$.
By first multiplying both sides of the Inequality \eqref{eqn:edge_bound} by
$(s-n\alpha_n)(n\alpha_1-s)/(\alpha_1-\alpha_n)^2$, we  algebraically manipulate this new
inequality to derive the following equivalent expression,
\begin{equation} \label{eqn:reordered_bound}
\frac{1}{4}(1+\alpha_1+\alpha_n)^2 - n\alpha_n \leq
\frac{\left( 2s - n(n-1) + (n-\alpha_n-1)(n-\alpha_n) -
\alpha_1(\alpha_1+1) \right)^2}{4(\alpha_1-\alpha_n)^2}.
\end{equation}

Combining the facts that $\alpha_1 \geq n-\alpha_n-1$ along with
$s \geq p\alpha_1 + (n-p)\alpha_n$,  this new Inequality
\eqref{eqn:reordered_bound} becomes
\begin{equation} \label{eqn:reordered_w_sum}
\begin{split}
\frac{1}{4}(1+\alpha_1+\alpha_n)^2 - n\alpha_n
&\leq
\frac{\left( 2s - n(n-1) + (n-\alpha_n-1)(n-\alpha_n) - \alpha_1(\alpha_1+1)
\right)^2}{4(\alpha_1-\alpha_n)^2} \\
&\leq
   \frac{\left( 2(p\alpha_1+(n-p)\alpha_n) - n(n-1) +
		   (n-\alpha_n-1)(n-\alpha_n) - \alpha_1(\alpha_1+1)
    \right)^2}{4(\alpha_1-\alpha_n)^2} \\
  &=
\left(p - \frac{1}{2}(1 + \alpha_1 + \alpha_n) \right)^2. \\
\end{split}
\end{equation}
We now rewrite the Inequality \eqref{eqn:reordered_w_sum} as 
\begin{equation} \label{eqn:reformulated}
p\alpha_1 \leq p(p-1) + (n-p)\alpha_n,
\end{equation}
implying $p\alpha_1 \leq  p(p-1) + \gamma + (n-p-1)\alpha_n$,
which is precisely the EG inequality for $\alpha'$ at $k=p$.

For $k=p+1$, we split the instance into two cases: $p \geq \gamma$
and $p < \gamma$.  If $p < \gamma$, then we use Theorem
\ref{thm:dur_index}, which states
that we only need to check to the largest index $k$ such that $\alpha_k \geq
k-1$, in order to prove that a sequence is graphic. 
Thus, if $p < \gamma$, then satisfying the EG
inequality for $k=p$ already established that the sequence is graphic.  Else,
if $p \geq \gamma$, then $\gamma \leq 2p - \alpha_{n-1}$ and summing
this inequality with the Inequality \eqref{eqn:reformulated}, we derive 
$p\alpha_1 + \gamma \leq (p+1)p + (n-(p+1))\alpha_n$, which is the 
EG inequality for $k=p+1$.
Therefore, the premise is established.
\end{proof}

\section{Remarks About The Bound}
This result generalizes the bound of Zverovich and Zverovich. We observe
that for the right hand side of Inequality \eqref{eqn:reordered_bound},
\begin{equation}
\frac{\left( 2s - n(n-1) + (n-\alpha_n-1)(n-\alpha_n) - \alpha_1(\alpha_1+1)
\right)^2}{4(\alpha_1-\alpha_n)^2} \geq 0.
\end{equation}
It follows that if $\frac{1}{4}(\alpha_1+\alpha_n+1)^2 - n\alpha_n \leq 0$ then the
sequence $\alpha$ is graphic by Theorem \ref{thm:edge_bound}.
By rearranging this inequality, we derive the bound of Theorem \ref{thm:zz}.

In general, this bound is sharp in the sense that there exist sequences that
satisfy the inequality, but any modification to parameters that causes the
inequality to be not satisfied allows for non-graphic sequences to be
formed.  For example, consider the following set of
degree sequences: for a given value of $\alpha_1$, define $\alpha_n = 2$,
$n = \alpha_1 +1$, and $s = 4\alpha_1-2$.  These values satisfy the
Inequality \eqref{eqn:edge_bound} with equality, and so any sequence with
these parameters is graphic. We notice
that modifying any one of these values so that the inequality is no
longer satisfied while keeping the other three fixed allows for non-graphic
sequences. For example,
if we increase $\alpha_1$ by 1 or decrease $n$ by 1, we have 
sequences, such as
$\seq{\alpha_1+1, \alpha_1-1, \underbrace{2,...,2}_{\alpha_1-1}}$ and
$\seq{\alpha_1, \alpha_1, 4,\underbrace{2,...,2}_{\alpha_1-3}}$ where 
the largest value in the sequence is larger than number of remaining values.
Likewise, if we decrease $\alpha_n$ by 1 or increase the sum by 2 (in order to
preserve an even sum), then we have sequences such as 
$\seq{\alpha_1, \alpha_1, 3, \underbrace{2,...,2}_{\alpha_1-3},1}$
and $\seq{\alpha_1, \alpha_1, 4,\underbrace{2,...,2}_{\alpha_1-2}}$ that again
can be quickly verified to be non-graphic.

In comparison to the Zverovich and Zverovich result,
there is a slight difference in the requirements for applying the bound.
This bound requires that $\alpha_1 \leq n-1$ and this restriction cannot be
relaxed.
This is because there are sequences where $\alpha_1>n-1$ but fulfill the
inequality (such as $\seq{8,6,6,6,6,6,6,6}$). In contrast to Zverovich
and Zverovich result, if that inequality is satisfied, it implies that
$\alpha_1 \leq n-1$ (proof of Theorem 1.2, \cite{Cairns:2016}).

It should be noted that this bound cannot be used for sets of regular sequences,
i.e., $\alpha_1 = \alpha_n$, but we can handle this case by using in
conjunction  the result that
states if $\alpha_1 - \alpha_n \leq 1$, then  $\alpha$ is graphic
(\cite{Chen:1988}, Lemma 1).  In addition, we note that when we extract
the smallest value in the sequence, we should choose
the smallest value where $\alpha_i > 0$.  This is because 
zeros in a degree sequence do not change whether the entire sequence is graphic
or not and can be discarded.

Finally, this result has practical use when
coupled with the usual Erd\H{o}s-Gallai conditions for
checking if a sequence is graphic. Before a sequence
can be tested using the EG conditions, it must
first be sorted. During this sorting step, the values for
largest and smallest indices, sequence length, and sequence sum can be easily
extracted. Thus, this new condition may save a number of steps in graphic
testing for certain sequences by verifying that a sequence is graphic without
having to compute the full EG inequalities (\cite{Cloteaux:2015}).

\section*{Acknowledgments}

The author would like to thank the anonymous referees for their
insightful comments and observations to improve the paper. Also, I would
like to thank Zolt\'an Szigeti for simplifying an argument that
ended up as Equation \eqref{eqn:p_inequality}.

\nocite{*}
\bibliographystyle{abbrvnat}
\bibliography{graphic_bounds}
\label{sec:biblio}

\end{document}